\documentclass[12pt]{article} 
\usepackage[utf8]{inputenc}
\usepackage{index}
\usepackage[reqno]{amsmath}
\usepackage{fixmath}
\usepackage{amssymb, amsthm, enumerate}
\usepackage{mathtools}
\usepackage{fullpage}
\usepackage{hyperref}
\usepackage[dvipsnames]{xcolor}

\usepackage{bbold}
\usepackage{authblk}
\usepackage{scalerel}
\newcommand\cprods{\mathbin{\scaleobj{0.8}{\square}}}
\usepackage{tikz}

\newcommand{\pa}[1]{\left(#1\right)} 
\newcommand{\IGNORE}[1]{}
\newcommand{\mat}[1]{\begin{matrix}#1\end{matrix}} 
\newcommand{\pmat}[1]{\pa{\mat{#1}}} 
\newtheoremstyle{plainsl}%
	{\topsep}
	{\topsep}
	{\slshape} 
	{}
	{\normalfont\bfseries}
	{.}
	{ }
	{}

\swapnumbers

{\theoremstyle{plainsl}
\newtheorem{theorem}{Theorem}[section]
\newtheorem{lemma}[theorem]{Lemma}
\newtheorem{prop}[theorem]{Proposition}
\newtheorem{corollary}[theorem]{Corollary}}
{\theoremstyle{remark}
}

\renewcommand\proof{\noindent\textsl{Proof. }}
\newcommand\sqr[2]{{\vbox{\hrule height.#2pt
    \hbox{\vrule width.#2pt height#1pt \kern#1pt
        \vrule width.#2pt}\hrule height.#2pt}}}

\renewcommand\qed{%
	\ifmmode\eqno\sqr53
	\else\nolinebreak\ \hfill\sqr53\medbreak\fi}

\title{Cayley graphs for extraspecial p-groups and a covering graph perspective on Huang's theorem}
\author{Maxwell Levit}
\affil{Department of Combinatorics \& Optimization, University of Waterloo}
\date{November 17, 2020}

\begin{document}
\maketitle

\begin{abstract}
In 1985, Arjeh Cohen and Jacques Tits proved the existence of a 4-cycle-free 2-fold cover of the hypercube. This Cohen-Tits cover is closely related to the signed adjacency matrix that Hao Huang used last year in his proof of the Sensitivity Conjecture. Terence Tao observed that Huang's signed adjacency matrix can be understood by lifting functions on an elementary abelian 2-group to functions on a central extension. Inspired by Tao's observation, we generalize the Cohen-Tits cover by constructing, as Cayley graphs for extraspecial $p$-groups, two infinite families of 4-cycle-free $p$-fold covers of the Cartesian product of $p$-cycles. 
\end{abstract}

\section{Introduction}
A key step in Hao Huang's proof of the Sensitivity Conjecture \cite{HH} is to construct a $\pm 1$ signing of the adjacency matrix of the hypercube in which the product of the signs around any 4-cycle is $-1$. The $d$-dimensional hypercube, denoted $Q_d$, is the Cartesian product of $d$ copies of the complete graph $K_2$, a graph that can reasonably be thought of as a degenerate 2-cycle. In this paper, for primes $p>2$, we will consider $p$th-root-of-unity valued ``signings'' of the adjacency matrices of Cartesian powers of $p$-cycles. In particular, we will generalize Huang's construction by producing labelings of these graphs in which the product of labels around any 4-cycle is not the identity.  

One of our main motivations is a covering graph perspective on signed graphs that we will now begin to describe. Let $X$ be a simple graph and $S$ a subset of its edges. The \textbf{signed graph} $(X,S)$ is the graph $X$ with edge labels $-1$ for $e\in S$ and $+1$ for $e\in E(X)\backslash S$. Given a signed graph $(X,S)$ we can construct a \textbf{2-fold cover} of $X$ by replacing each vertex $v$ of $X$ with two vertices, $v'$ and $v''$ and replacing each edge $\{u,v\}$ with a pair of edges: $\{\{v'u''\},\{v''u'\}\}$ if $\{u,v\}\in S$, or $\{\{v'u'\},\{v''u''\}\}$ if $\{u,v\}\notin S$. The choice of labels $\pm1$ is not cosmetic: Let $A$ denote the adjacency matrix of $X$, and $\hat{A}$ the \textbf{signed adjacency matrix}

\[\hat{A}_{i,j}=
\begin{cases} ~~0 &\text{if}~~~ ij\notin E(X)\\  
~~1 &\text{if}~~~ ij\in E(X)\backslash S\\
 -1 &\text{if}~~~ ij\in S
\end{cases}.\]

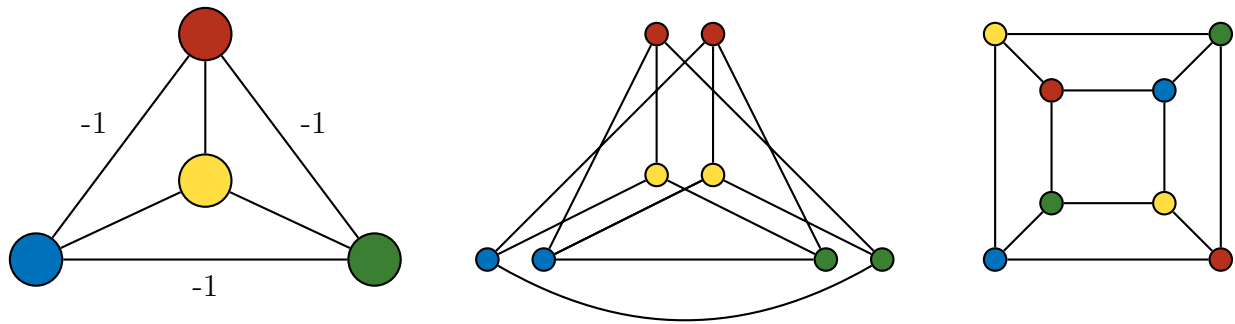
\begin{figure}
    \centering
\tikzstyle{every node}=[circle, draw, fill=white!50,
                        inner sep=7pt, minimum width=4pt]
\def\x{8.5}
\begin{tikzpicture}[thick,scale=.75]%
\node (a) at (0,0) [fill=RoyalBlue] {};
\node (b) at (3,4)[fill=BrickRed]{};
\node (c) at (6,0)[fill=OliveGreen]{};
\node (d) at (3,1.4)[fill=Goldenrod]{};

\tikzstyle{every node}=[circle,
                        inner sep=9pt, minimum width=4pt]
\draw{ 
    (a)--(b) node [midway, xshift=-10pt, yshift=9pt] {-1}
    (a)--(c) node [midway, yshift=-10pt] {-1}
    (a)--(d)
    (b)--(c) node [midway, xshift=9pt, yshift=9pt] {-1}
    (b)--(d) 
    (c)--(d) 
};
   
\tikzstyle{every node}=[circle, draw,
                        inner sep=3pt, minimum width=4pt]
\node (1) at (-.5+\x,0)[fill=RoyalBlue] {};
\node (2) at (.5+\x,0)[fill=RoyalBlue] {};
\node (3) at (3.5+\x,4)[fill=BrickRed] {};
\node (4) at (2.5+\x,4)[fill=BrickRed] {};
\node (5) at (2.5+\x,1.5) [fill=Goldenrod]{};
\node (6) at (3.5+\x,1.5) [fill=Goldenrod]{};
\node (7) at (5.5+\x,0)[fill=OliveGreen] {};
\node (8) at (6.5+\x,0)[fill=OliveGreen] {};
\draw{  
        (1)--(3)
        (1) edge[bend right](8)
        (1)--(5)
        (2)--(7)
        (2)--(4)
        (2)--(6)
        (3)--(7)
        (4)--(8)
        (4)--(5)
        (5)--(7)
        (6)--(3)
        (6)--(8)

};
\draw{(2)--(6)};

\node (q1) at (17,0)[fill=RoyalBlue] {};
\node (q2) at (21,0) [fill=BrickRed] {}; 
\node (q3) at (17,4)[fill=Goldenrod] {}; 
\node (q4) at (21,4)[fill=OliveGreen] {}; 
 \node (q5) at (18,1)[fill=OliveGreen] {}; 
 \node (q6) at (18,3)[fill=BrickRed] {}; 
 \node (q7) at (20,1)[fill=Goldenrod] {}; 
 \node (q8) at (20,3)[fill=RoyalBlue] {}; 
 
 \draw{ 
 (q1)--(q2)
 (q1)--(q3)
 (q2)--(q4)
 (q3)--(q4)
 
 (q5)--(q6)
 (q5)--(q7)
 (q6)--(q8)
 (q7)--(q8)
 
 (q1)--(q5)
 (q2)--(q7)
 (q3)--(q6)
 (q4)--(q8)};

\end{tikzpicture}

    \caption{A signed $K_4$  (with $+1$ labels omitted) and the associated 2-fold cover, which is isomorphic to the cube.}
\end{figure}

One can verify that the (unsigned) adjacency matrix of the associated 2-fold cover of $X$ is similar to the direct sum $A \oplus \hat{A}$. This demonstrates a close connection between 2-fold covering graphs and the spectra of $\pm1$ signings of adjacency matrices. Probably the most notable application of this connection is the venture proposed by Bilu and Linial \cite{BL}: To construct expander graphs from an iterated sequence of 2-fold covers. Marcus, Spielman, and Srivastava \cite{MSS} have seen this venture to completion with resounding success.

Another recent triumph of the signed adjacency matrix can be found in Huang's remarkable proof \cite{HH} of the Sensitivity Conjecture. This conjecture of Nisan and Szegedy \cite{NS} asked if two measures of the complexity of a Boolean function are polynomially related. It was translated by Gotsman and Linial \cite{GL} to a problem of lower bounding the maximum degree of an induced subgraph of the hypercube, $Q_d$. Huang solved this version of the problem, proving the Sensitivity Conjecture. A key step in his solution is to construct a very particular $\pm1$ signing of the adjacency matrix of $Q_d$.

As it turns out, the 2-fold cover of $Q_d$ associated with Huang's matrix is a remarkable graph in its own right. This graph was considered, under very different circumstances, by Arjeh Cohen and Jacques Tits \cite{CT} in 1985. Cohen and Tits were interested in proving the uniqueness of several notable incidence systems: The generalized hexagons of orders (2,2) and (2,8), and the near octagon of order (2,4;0,3). Their investigation into the girth six 2-fold cover of $Q_d$ was prompted by its presence as a subgraph of the colinearity graph of these incidence systems.
We remark that this is not the only context within finite geometry where these graphs appear: The Cohen-Tits cover of $Q_3$ is isomorphic to the point-line incidence graph of the M\"obius-Kantor configuration. 

In Section 2 we will verify that the 2-fold cover associated with Huang's signed adjacency matrix of $Q_d$ the Cohen-Tits cover. This connection was observed by Godsil, Silina and the author in \cite{GLS}. In Section 3 we then present a deeper connection between these objects. We show that the Cohen-Tits cover is a Cayley graph for a central extension of $\mathbb{Z}_2^d$ by $\mathbb{Z}_2$. In particular, this central extension is the Heisenberg extension of $\mathbb{Z}_2^d$ which Tao uses in \cite{TT} to describe the twisted convolution approach to Huang's theorem. This will lead to our main result, a generalization of the Cohen-Tits covers to Cartesian products of cycles of prime length.

\begin{theorem}\label{MT}
Let $p$ be an odd prime. Let $X$ be a Cartesian product of $d>1$ cycles of length $p$. There are (at least) two non-isomorphic $p$-fold covers of $X$ containing no 4-cycles.
\end{theorem}

We will construct these covers as Cayley graphs for central extensions of $\mathbb{Z}_p^d$, namely the extraspecial $p$-groups. This allows us to recast the combinatorial condition of being 4-cycle-free as the non-commutiativity of a subset of group elements. Moreover, this allows us  (in Section 6) to derive the promised $p$th-root-of-unity valued signing of the Cartesian products from the 2-cocycle defining this central extension. We believe a similar approach will be useful in constructing other interesting covers in the future. 

We probably owe the reader some indication of what makes our newly constructed covers interesting. One answer is that, in some cases, their adjacency matrices sit inside coherent algebras or association schemes. In other cases these covers have as few automorphisms as possible (given that they are Cayley graphs), which is nontrivial, given the large amount of symmetry in the graphs they cover. We will say more about this in Section 7.

Also in Section 7 we discuss the possibility of extracting useful spectral information from our covers. These covers, or the associated $p$th-root-of-unity valued signings, are analogous to Huang's matrix. So one might hope that an analog of Huang's proof will provide good bounds on the maximum degree of induced subgraphs of Cartesian products of $p$-cycles. Indeed, for small examples we are able to improve slightly upon a result of Tikaradze \cite{AT}. However we expect that, in general, our covers do not have large enough eigenvalues of high enough multiplicity to yield useful bounds.

\section{Covering graphs}

The informal definition of 2-fold covering graphs given in the introduction can be extended to a broader notion of covering graphs as follows. Let $Y=(V,E)$ be a (simple) graph. The graph $X$ is a \textbf{cover of} $Y$ if there is a graph homomorphism $\gamma:X\rightarrow Y$ (called the \textbf{covering map}), satisfying

\begin{enumerate}[a)]
    \item Each \textbf{fiber} $\gamma^{-1}(v)$ is an independent set in $X$. 
    \item For each edge $uv\in E(Y)$, the subgraph of $X$ induced by $\gamma^{-1}(v)\cup \gamma^{-1}(u)$ is a perfect matching.

\end{enumerate}

If $Y$ is connected, all fibers must have the same size, say $r$. In this case we say that $X$ is an $r$-\textbf{fold} cover. 

Our main result is a generalization of the following result of Cohen and Tits.

\begin{theorem}[\cite{CT} Lemma 4a]
There exists a 4-cycle-free, 2-fold cover of the hypercube $Q_d$.
\end{theorem}

Cohen and Tits prove this result by analyzing a particular subgroup of index 2 within the fundamental group of the (1-skeleton of) the hypercube. The wonderful correspondence between subgroups of the fundamental group and covering spaces then yields the result. See \cite{AH} Section 1.3 for a treatment of this correspondence in general, and \cite{CT} for this particular application. We remark on this aspect of their proof because it is quite elegant and very natural from a topological perspective. For our purposes it is instructive to give an elementary inductive proof.

\proof{We wish to construct a 2-fold cover of $Q_d$ which contains no 4-cycles. On the merits of our discussion from Section 1, it suffices to construct a signed adjacency matrix for $Q_d$ in which every 4-cycle contains an odd number of edges with label $-1$. This ensures that the associated 2-fold cover contains no 4-cycles. To find such a signing, we employ the inductive construction of $Q_d$ as two copies of $Q_{d-1}$ joined by a perfect matching. 
 
Suppose we have found a signed adjacency matrix $A_{d-1}$ for $Q_{d-1}$ with the desired property. Notice that the ``opposite'' signing, $-A_{d-1}$, in which all $1$'s have been changed to $-1$'s and vice versa, will also have the desired property. Joining two oppositely signed copies of $Q_{d-1}$ by a perfect matching whose edges are all given the same sign (say $+1$) yields a signing of $A(Q_d)$ in which each 4-cycle has an odd number of $-1$ edges. This in turn yields a 4-cycle-free cover of $Q_d$.} \qed

We obtain the connection to Huang's matrix as an immediate corollary to this proof.
\begin{corollary}
Let $A_d$ be the signed adjacency matrix of $Q_d$ obtained in Lemma 2 of \cite{HH}. The 2-fold cover of $Q_d$ associated with $A_d$ is the Cohen-Tits cover.
\end{corollary}

\proof{When we translate ``joining two oppositely signed copies of $Q_{d-1}$ by a perfect matching'' to operations on the associated signed adjacency matrices, we obtain exactly \[A_1=\pmat{0&1\\1&0},~~~A_{d}=\pmat{A_{d-1}&I_{d-1}\\I_{d-1}&-A_{d-1}}.\]}\qed

Cohen and Tits's result also shows that the 4-cycle free cover of $Q_d$ is unique up to isomorphism. This result will be useful for us later.

\begin{theorem}[\cite{CT} Lemma 4a]
Up to isomorphism there is only one 4-cycle free 2-fold cover of $Q_d$. \qed
\end{theorem}

\section{Convolutions and Cayley graphs}

Let $\mathbb{Z}_p^d$ denote the additive group of a $d$-dimensional vector space over the field with $p$ elements. In this section we show that the Cohen-Tits cover of $Q_d$ is a Cayley graph for $H_d$, the $\textbf{Heisenberg extension}$ of $\mathbb{Z}_2^d$ that is defined as follows. Let $\beta:\mathbb{Z}_2^d\times \mathbb{Z}_2^d\rightarrow \mathbb{Z}_2$  be the bilinear form given by  \[\beta(x,y):=\sum_{1\leq i <j \leq d} x_iy_j.\] The group $H_d$ is defined on the set $\mathbb{Z}_2^d\times \mathbb{Z}_2$ with multiplication \[(x,t)(y,s)=(x+y,s+t+\beta(x,y)).\]

Terence Tao observed the connection between $H_d$ and Huang's signed adjacency matrix in his blog post \cite{TT}. This connection inspired Theorem 3.1 of this section, which in turn inspired the proofs of our main results.  We now quickly summarize the relevant details of Tao's exposition.

Let $G$ be a finite group, and let $\ell^2(G)$ be the space of functions from $G$ to $\mathbb{C}$. For $f,g\in \ell^2(G)$, the \textbf{convolution} $\star_G: \ell^2(G)\times \ell^2(G)\rightarrow \ell^2(G)$ is the function defined by \[f\star_G g(x)=\sum_{y\in G}f(y)g(y^{-1}x).\]

Tao interprets the action of the adjacency matrix of $Q_d$ on $f\in \ell^2(\mathbb{Z}_2^d)$ as the convolution  $Af=f\star_{\mathbb{Z}_2^d} \mu$, where $\mu(x)=1$ if $x$ is any of the standard basis vectors $e_1,\dots e_d$, and $\mu(x)=0$ otherwise.  He then defines the \textbf{twisted convolution} $\star_\beta: \ell^2(\mathbb{Z}_2^d)\times \ell^2(\mathbb{Z}_2^d)\rightarrow \ell^2(\mathbb{Z}_2^d)$ given by \[f\star_\beta g(x)=\sum_{y\in \mathbb{Z}_2^d}(-1)^{\beta(y,y^{-1}x)}f(y)g(y^{-1}x)\] and an operator $A_\beta$ whose action on $f$ is defined by $A_\beta f=f\star_\beta \mu$. The operator $A_\beta$ has the same spectrum as Huang's signed adjacency matrix and can be used to derive Huang's result. 
 
Tao then remarks upon the map $L:\ell^2(\mathbb{Z}_2^d)\rightarrow \ell^2(H_d)$, defined by $f\mapsto L(f)$ with  \[ L(f)(x,t)=(-1)^t f(x).\] One can show that for $f,g\in \ell^2(\mathbb{Z}_2^d)$,\[L(f)\star_{H_d} L(g)=L(f\star_\beta g)\] which relates Huang's matrix to $H_d$ via the lift $L$. Theorem 3.1  exposes this relationship from a slightly different perspective.

Given a group $G$ with identity $e$ and an inverse closed subset $S\subseteq G\backslash \{e\}$, the \textbf{Cayley graph} $\text{Cay}(G,S)$ is the graph with vertex set $G$ and $\{g,h\}$ an edge if $gh^{-1}\in S$. One very pertinent example is the hypercube, which is isomorphic to the Cayley graph $\text{Cay}(\mathbb{Z}_2^d,S)$ with $S=\{e_1,\dots, e_d\}$, the standard basis for $\mathbb{Z}_2^d$.

\begin{theorem}
Let $S_d=\{(e_1,0),(e_2,0)\dots,(e_d,0)\}$. The Cohen-Tits cover of $Q_d$ is isomorphic to $\text{Cay}(H_d,S_d)$.
\end{theorem}

\proof{First note that the connection set $S_d$ is inverse closed and does not contain the identity, so our Cayley graph is well defined. The projection $\gamma:H_d\rightarrow \mathbb{Z}_2^d$ defined by $\gamma(g,s)=g$ serves as a 2-fold covering map from $\text{Cay}(H_d,S_d)$ onto $\text{Cay}(\mathbb{Z}_2^d,\gamma(S_d))\cong Q_d$. In deference to Theorem 2.3, all that remains is to show that $\text{Cay}(H_d,S_d)$  has no 4-cycles.

Suppose $\text{Cay}(H_d,S_d)$ contains a 4-cycle $C$. This cycle's image $\gamma(C)$ must be a 4-cycle in $\text{Cay}(\mathbb{Z}_2^d,\gamma(S_d))$. So there is some $g\in \mathbb{Z}_2^d$ and $i,j\in \{1,\dots n\}$, with $i<j$, so that \[V(\gamma(C))=\{g,g+e_i,g+e_i+e_j,g+e_j\}.\]

This implies the existence of $t\in 
\{0,1\}$ so that \[V(C)=\{(g,t),(g,t)(e_i,0),(g,t)(e_i,0)(e_j,0),(g,t)(e_i,0)(e_j,0)(e_i,0)^{-1}\}.\]

But since $C$ is a 4-cycle we have \[(e_i,0)(e_j,0)(e_i,0)^{-1}(e_j,0)^{-1}=(\bar{0},0).\]

This contradicts the following direct calculation from the definition of $H_d$. 
\begin{align*}
(e_i,0)(e_j,0)(e_i,0)^{-1}(e_j,0)^{-1}
&=(e_i+e_j,\beta(e_i,e_j))(e_i,0)(e_j,0)\\
&=(e_i,\beta(e_i,e_j)+\beta(e_i+e_j,e_i))(e_i,0)\\
&=(\bar{0},\beta(e_i,e_j)+\beta(e_i+e_j,e_i)+\beta(e_i,e_i)\\
&=(\bar{0},1+0+0)\\
&=(\bar{0},1).\\
\end{align*}}\qed

\noindent\textbf{Remark:} One can construct an eigenbasis for the adjacency matrix of the Cohen-Tits cover by taking a full collection of the eigenfunctions for $A(Q_d)$ and ``doubling them up'' into functions which are constant on the fibers of the cover. This accounts for half the eigenvaules of the cover. The other half belong to eigenspaces for the ``new'' eigenvalues $\sqrt{d}$ and $-\sqrt{d}$. These new eigenfunctions are orthogonal to the old ones, so for any such function $f$ and fiber $\{v,v'\}$ we have $f(v)=-f(v')$. It follows that the map $L$ lifts any function on $Q_d$ to a function on the cover within the support of the $\sqrt{d}$ and $-\sqrt{d}$ eigenspaces.

\section{Extraspecial $p$-groups}

We have identified the Cohen-Tits covers as Cayley graphs for a central extension of $\mathbb{Z}_2^d$. A natural next step is to consider Cayley graphs for central extensions of $\mathbb{Z}_p^d$, $p$ an odd prime. To this end, we consider the extraspecial $p$-groups, whose relevant properties we now recall. Basic definitions and more details can be found in many places, e.g. \cite{DG}, \cite{MS}.

Fix a prime $p$. A $p$-group $G$ is \textbf{extraspecial} if its center $Z=Z(G)$ is a cyclic group of order $p$ and $G/Z$ is a (non-trivial) elementary abelian $p$-group. For each nonnegative integer $d$ there are exactly two non-isomorphic extraspecial $p$-groups of order $p^{1+2d}$ and none of order $p^{2d}$ (see \cite{DG} Section 5.2). One of these groups has exponent $p$, the other has exponent $p^2$. We follow the convention to denote these groups by $p_+^{1+2d}$ and $p_-^{1+2d}$, respectively.

 It follows from the definition that the extraspecial $p$-groups admit short exact sequences \[1\rightarrow \mathbb{Z}_p\rightarrow p_\pm^{1+2d}\rightarrow \mathbb{Z}_p^{2d}\rightarrow 1\] in which $\mathbb{Z}_p$ is identified with $Z=Z(p_\pm^{1+2d})$. A short exact sequence of this type is known as \textbf{central extension} of $\mathbb{Z}_p^{2d}$ by $\mathbb{Z}_p$, and can be specified by defining a \textbf{2-cocycle}, a map $\kappa: \mathbb{Z}_p^{2d}\times \mathbb{Z}_p^{2d}\rightarrow Z$ satisfying \[\kappa(a+b,c)+\kappa(a,b)=\kappa(a,b+c)+\kappa(b,c).\] This 2-cocycle $\kappa$ can be used to explicitly define the third group in the short exact sequence above via the following multiplication on $\mathbb{Z}_p^{2d}\times Z$. \[(g,s)(h,t)=(g+h,s+t+\kappa(g,h)).\] The next two lemmas explicate 2-cocycles that give rise to the two isomorphism classes of extraspecial $p$-groups. We will identify $\mathbb{Z}_p^{2d}$ with $\mathbb{Z}_p^d\times\mathbb{Z}_p^d$ in order to make the notation a bit cleaner.

\begin{lemma}

The map $\kappa_+: (\mathbb{Z}_p^{d}\times \mathbb{Z}_p^{d})\times (\mathbb{Z}_p^{d}\times \mathbb{Z}_p^{d})\rightarrow \mathbb{Z}_p$ defined by \[((\bar{a},\bar{b}),(\bar{c},\bar{d}))\mapsto \bar{b}\cdot \bar{c}\] is a 2-cocycle. It determines the extraspecial $p$-group $p_+^{1+2d}$ of exponent $p$ on the set  $\mathbb{Z}_p^{d}\times \mathbb{Z}_p^d\times  \mathbb{Z}_p$.
\end{lemma} 

\proof{The cocycle condition follows from the linearity of the dot product. The induced group multiplication is

 \begin{equation}\label{plusmult}
     (\bar{a},\bar{b},z)(\bar{c},\bar{d},w)=(\bar{a}+\bar{c},\bar{b}+\bar{d},z+w+\bar{b}\cdot \bar{c}).
 \end{equation}

Associativity of this multiplication follows from the cocycle condition. The identity element is $(\bar{0},\bar{0},0)$ and the inverse is given by \[(\bar{a},\bar{b},z)^{-1}=(-\bar{a},-\bar{b},-z+\bar{a}\cdot \bar{b}).\]

Denote this group by $G$. To verify that $G$ is extraspecial, we calculate the commutator of a generic pair of elements.

\begin{align*}
    [(\bar{a},\bar{b},z),(\bar{c},\bar{d},w)]
    &=(\bar{a},\bar{b},z)^{-1}(\bar{c},\bar{d},w)^{-1}(\bar{a},\bar{b},z)(\bar{c},\bar{d},w)\\
    &=(-\bar{a},-\bar{b},-z+\bar{a}\cdot \bar{b})(-\bar{c},-\bar{d},-w+\bar{c}\cdot \bar{d})(\bar{a},\bar{b},z)(\bar{c},\bar{d},w)\\
    &=(0,0,\bar{a}\cdot\bar{b}+\bar{c}\cdot \bar{d}+(-\bar{b})\cdot (-\bar{c})+(-\bar{b}-\bar{d})\cdot\bar{a}+(-\bar{d})\cdot\bar{c})\\
    &=(0,0,\bar{b}\cdot \bar{c}-\bar{a}\cdot\bar{d}).
\end{align*}

It follows that a fixed element $(\bar{a},\bar{b},z)$ is in the center $Z=Z(G)$ if and only if $\bar{b}\cdot \bar{x}-\bar{y}\cdot \bar{a}$ is identically zero as $(\bar{x},\bar{y})$ ranges over $\mathbb{Z}_p^d\times \mathbb{Z}_p^d$. This happens if and only if $\bar{a}=\bar{b}=\bar{0}$, hence
 \[Z=\{(\bar{0},\bar{0},z):z\in \mathbb{Z}_p\}.\] Next, note that the cosets $gZ$ in $G$ depend only on the first $2d$ coordinates of the vector $g$. After we restrict to these coordinates the multiplication rule for the group degenerates to vector addition, so $G/Z$ is isomorphic to $\mathbb{Z}_p^{2d}$.
 
Finally, we calculate \[(\bar{a},\bar{b},x)^p=(\bar{0},\bar{0},\bar{b}\cdot \bar{a}+2\bar{b}\cdot a+\dots+(p-1)\bar{b}\cdot \bar{a})= (\bar{0},\bar{0},\binom{p}{2} \bar{b}\cdot \bar{a})=(\bar{0},\bar{0},0).\] So the exponent is $p$ and $G$ is $p_+^{1+2d}$.}\qed

Next we give a similar lemma describing the other extraspecial group of order $p^{1+2d}$ over the same set. This requires a less natural multiplication. Let $\iota$ be the inclusion $\mathbb{Z}_p\xhookrightarrow{} \mathbb{Z}$. Define a map  $\phi:\mathbb{Z}_p\times \mathbb{Z}_p \rightarrow \mathbb{Z}_p$ by \[\phi(a,b)=\begin{cases} 0, &\text{if } \iota(a)+\iota(b)<p \\ 1, &\text{if } \iota(a)+\iota(b)\geq p. \end{cases}\] 

Let $\bar{a}, \bar{c}$ be vectors in $\mathbb{Z}_p^d$, and let $a_1$, $c_1$ denote their first coordinates. We use the function $\phi$ to ``carry'' the data from this first coordinate to the coordinate of the commutator.

\begin{lemma}
The map $\kappa_-: (\mathbb{Z}_p^{d}\times \mathbb{Z}_p^{d})\times (\mathbb{Z}_p^{d}\times \mathbb{Z}_p^{d})\rightarrow \mathbb{Z}_p$ defined by \[((\bar{a},\bar{b}),(\bar{c},\bar{d}))\mapsto \bar{b}\cdot \bar{c}+\phi(a_1,c_1)\] is a 2-cocycle. It determines the extraspecial $p$-group, $p_-^{1+2d}$ of exponent $p^2$ on the set  $\mathbb{Z}_p^{d}\times \mathbb{Z}_p^d\times  \mathbb{Z}_p$.
\end{lemma}

\proof{As in the previous lemma, the first summand in the expression for $\kappa_-$ will satisfy the cocycle condition since the dot product is linear. For the second summand we must show that 

\[\phi(a+b,c)+\phi(a,b)=\phi(a,b+c)+\phi(b,c).\]

If $\iota(a)+\iota(b)+\iota(c)<p$ then all four terms in this equation are 0. If $2p<\iota(a)+\iota(b)+\iota(c)$ then all four terms are 1. Now suppose $p<\iota(a)+\iota(b)+\iota(c)<2p$. If $\iota(a)+\iota(b)<p$ then $\iota(a+b)+\iota(c)\geq p$, and the left side of the equation is $1+0$. If $\iota(a)+\iota(b)\geq p$ then $\iota(a+b)+\iota(c)<p$, and the left side of the equation is $0+1$. Hence the left side is identically 1. The same argument applied to $\iota(b)+\iota(c)$ shows that the right side is also identically 1. So $\kappa_-$ is a 2-cocycle. The induced multiplication is 
\begin{equation}\label{minusmult}
    (\bar{a},\bar{b},z)(\bar{c},\bar{d},w)=(\bar{a}+\bar{c},\bar{b}+\bar{d},z+w+\bar{b}\cdot \bar{c}+\phi(a_1,c_1)).
\end{equation}

Associativity of this multiplication follows from the cocycle condition. The identity element is $(\bar{0},\bar{0},0)$ and the inverse is given by \[(\bar{a},\bar{b},z)^{-1}=(-\bar{a},-\bar{b},-z+\bar{a}\cdot \bar{b}-\phi(a_1,-a_1)).\]

Denote this group by $G$. To verify that $G$ is extraspecial, we calculate the commutator of a generic pair of elements.
\begin{align*}
    &[(\bar{a},\bar{b},z),(\bar{c},\bar{d},z)]\\&~=(\bar{a},\bar{b},z)^{-1}(\bar{c},\bar{d},z)^{-1}(\bar{a},\bar{b},z)(\bar{c},\bar{d},z)\\
    &~=(-\bar{a},-\bar{b},-z+\bar{a}\cdot \bar{b}-\phi(-a_1,-a_1))(-\bar{c},-\bar{d},-z'+\bar{c}\cdot \bar{d}-\phi(c_1,-c_1))(\bar{a},\bar{b},z)(\bar{c},\bar{d},z)\\
    &~=(0,0,\bar{b}\cdot \bar{c}-\bar{a}\cdot\bar{d}-\phi(c_1,-c_1)-\phi(-a_1,a_1)+\phi(-a_1,-c_1)+\phi(-a_1-c_1,a_1))+\phi(-c_1,c_1)\\
    &~=(0,0,\bar{b}\cdot \bar{c}-\bar{a}\cdot\bar{d}-\phi(-a_1,a_1)+\phi(-a_1,-c_1)+\phi(-a_1-c_1,a_1)). 
\end{align*}

We claim that \[-\phi(-a_1,a_1) +\phi(-a_1,-c_1)+\phi(-a_1-c_1,a_1)=0\] for all values of $a_1$ and $c_1$.  This is immediate when $a_1=0$ since all three summands are 0. When $a_1>0$, we have $-\phi(-a_1,a_1)=-1$ and it remains to verify that $\phi(-a_1,-c_1)$, and $\phi(-a_1-c_1,a_1)$ take different values.

Suppose $\phi(-a_1,-c_1)=1$. Then $\iota(-a)+\iota(-c)\geq p$ which implies $\iota(-a-c)=\iota(-a)+\iota(-c)-\iota(p)$, hence $\iota(-a-c)+\iota(a)=\iota(c)+\iota(-a)+\iota(a)-\iota(p)=\iota(c)<p$, and $\phi(-a-c,a)=0$.

Suppose $\phi(-a_1,-c_1)=0$. Then $\iota(-a)+\iota(-c)< p$  which implies $\iota(-a-c)=\iota(-a)+\iota(-c)$, hence $\iota(-a-c)+\iota(a)=\iota(c)+\iota(-a)+\iota(a)=p+\iota(c)\geq p$, and $\phi(-a-c,a)=1$. This proves the claim, and reduces our commutator to the familiar formula

\begin{equation*}
    [(\bar{a},\bar{b},z),(\bar{c},\bar{d},z)]=(0,0,\bar{b}\cdot \bar{c}-\bar{a}\cdot\bar{d}).
\end{equation*}

As in the proof of Lemma 4.1, this implies that $Z(G)=\{(\bar{0},\bar{0},z):z\in \mathbb{Z}_p\}$. This in turn implies that $G/Z$ is elementary abelian. 

Finally, let $\bar{v}$ be the vector $(1,0,\dots 0)$ in $\mathbb{Z}_p^d$, we calculate $(\bar{v},\bar{0},0)^p$.
Note that for $k<p$, $\phi(k-1,1)=0$, and so  $(\bar{v},\bar{0},0)^k=(k\bar{v},\bar{0},0)$. Hence \[(\bar{v},\bar{0},0)^p=((p-1)\bar{v},\bar{0},0)(\bar{v},\bar{0},0)=(\bar{0},\bar{0},1).\] So $(\bar{v},\bar{0},0)$ has order $p^2$. It follows that $G$ has exponent $p^2$, and is $p_-^{1+2d}$.} \qed 

\noindent\textbf{Remark:} As a point of interest with no real bearing on what follows, we recall a definition of Philip Hall \cite{PH}. Let $[G,G]$ denote the group generated by all commutators of $G$. Two groups $G_1$ and $G_2$ are \textbf{isoclinic} if there are isomorphisms $\theta:G_1/Z(G_1)\rightarrow G_2/Z(G_2)$ and $\psi:[G_1,G_1]\rightarrow [G_2,G_2]$ so that for all $g,h\in G_1/Z(G_1)$, \[ \psi([g,h])=[\theta(g),\theta(h)].\]

By expressing  $p_+^{1+2d}$ and $p_-^{1+2d}$ as groups on the same set with identical commutators, we have shown that these groups are isoclinic by choosing $\theta$ and $\psi$ to be identity maps.

\section{Our construction}

Given graphs $X=(V_1,E_1), Y=(V_2,E_2)$, the \textbf{Cartesian product} $X\cprods{}Y$ is the graph with vertex set $V_1\times V_2$ and $(v_1,v_2)$ adjacent to $(u_1,u_2)$ if either \[v_1=u_1 \text{ and } v_2\sim_Yu_2\] or \[v_2=u_2 \text{ and } v_1\sim_Xu_1.\] 

Let $X^d$ denote the Cartesian product of $d$ copies of $X$. Let $C_p$ denote the cycle with $p$ vertices. Fix an odd prime $p$ and a non-negative integer $d$. In this section we construct two non-isomorphic 4-cycle-free $p$-fold covers of $C_p^{2d}$. 

Let $\{e_1,\dots,e_{d}\}$ be the standard basis for $\mathbb{Z}_p^{d}$. It is well known (e.g. \cite{WLG}) that $C_p^{d}$ is a Cayley graph for $\mathbb{Z}_p^{d}$ with connection set $\{\pm e_1,\dots, \pm e_d\}$. We will employ this description of $C_p^d$ along with the theory of extraspecial $p$-groups in order to prove Theorem \ref{MT}. In particular, we will construct our covers as Cayley graphs for the extraspecial $p$-groups discussed in Section 4. This will amount to specifying connection sets $S_\pm$ with several nice properties. 

Let $e_1,\dots e_d,f_1,\dots f_d$ denote the standard basis for $\mathbb{Z}_p^d\times \mathbb{Z}_p^d$, and let $[n]=\{1,\dots n\}$.  Define sets
 \[A=\left\{ \sum_{i\in [k]} e_i +\sum_{j\in [k-1]} f_j ~:~   k \in [d]\right\},~~~ B=\left\{ \sum_{i\in [k-1]} e_i +2e_{k}+\sum_{j\in [k]} f_j  ~:~  k \in [d]\right\}\]
 
and \[S=A\cup B.\]
 
For example, if $d=2$, $S$ consists of:
\begin{center}
$(1,0)\times(0,0)\in A$,

$(2,0)\times(1,0)\in B$,

$(1,1)\times(1,0)\in A$,

$(1,2)\times(1,1)\in B$.
\end{center}

Recall from the previous section that we may consider $\mathbb{Z}_p^d\times \mathbb{Z}_p^d\times \mathbb{Z}_p$ as the set on which both extraspecial $p$-groups $p_\pm^{1+2d}$ are defined. We will use the set $S$ to build a connection set for Cayley graphs on these groups. To this end, let $\epsilon: \mathbb{Z}_p^d\times \mathbb{Z}_p^d\xhookrightarrow{}  \mathbb{Z}_p^d\times \mathbb{Z}_p^d\times \mathbb{Z}_p$ be the natural inclusion defined by $\epsilon(\hat{a},\hat{b})=(\hat{a},\hat{b},0)$.

The next lemma will ensure that the set $\epsilon(S)$ has the necessary properties to produce covers of $C_p^{2d}$ that contain no 4-cycles.

\begin{lemma}
Let $S=A\cup B$ and $\epsilon$ be defined as above.
\begin{enumerate}[a)]
\item $S$ is a basis for the vector space  $\mathbb{Z}_p^{2d}$.
\item No two elements of $\epsilon(S)$ commute with respect to the multiplication (\ref{plusmult}) or (\ref{minusmult}).

\end{enumerate}
\end{lemma}

\proof{Write the elements of $S$ as the rows of a $2d\times 2d$ matrix. Move all the even indexed columns to the end while preserving the relative order of these columns. In other words, apply the permutation of column indices given by \[2i\mapsto d+i, ~2d-i\mapsto d-i\] for $i\in \{1,\dots, d\}$. It is immediate that this matrix is lower triangular with non-zero diagonal entries, hence its rows are a basis for $\mathbb{Z}_p^{2d}$. This proves a).

We have noted in the previous section that the commutator maps for $p_-^{1+2d}$ and $p_+^{1+2d}$ are identical. To prove b) it suffices to verify that this commutator is nonzero for distinct $g,h\in \epsilon(S)$. Indeed, we calculate

\begin{center}
$[g,h]=\begin{cases} 
(0,0,1), & g,h\in A\\ (0,0,-1), & g,h\in B \\(0,0,1), &g\in A,~h\in B.
\end{cases}$
\end{center}
}~\qed

The final step in our construction is to make $\epsilon(S)$ into a connection set for a Cayley graph of $p_\pm^{1+2d}$. To do this we must take the set's closure under inverse with respect to the appropriate multiplication. Define $S_+$, respectively $S_-$, by \[S_\pm=\epsilon(S)\cup \epsilon(S)^{-1}\] where $\epsilon(S)^{-1}$ consists of the inverses of elements of $\epsilon(S)$ with respect to the multiplication (\ref{plusmult}), respectively (\ref{minusmult}).

\begin{theorem}
 Let $p$ be an odd prime and $d$ a non-negative integer. Let $S$, $\epsilon$, $S_\pm$ be defined as above. The Cayley graphs $\text{Cay}(p_+^{1+2d},S_+)$ and $\text{Cay}(p_-^{1+2d},S_-)$ are non-isomorphic 4-cycle-free $p$-fold covers of $C_p^{2d}$.
\end{theorem}

\proof{Let $T=\{e_1,\dots,e_{2d}\}$ be the standard basis for $\mathbb{Z}_p^{2d}$. Recall that $C_p^{2d}$ is a Cayley graph for $\mathbb{Z}_p^{2d}$ with connection set \[T\cup T^{-1}=\{\pm e_1,\dots, \pm e_d\}.\] $T$ is a basis for $\mathbb{Z}_p^{2d}$ and by Lemma 5.1 a) the set $S$ is as well. Choose a bijection from $T$ to $S$, and let $\alpha$ be its extension to a linear map. It follows that $\alpha$, considered as a map between vertex sets, is an isomorphism from $\text{Cay}(\mathbb{Z}_p^{2d},T\cup T^{-1})$ to $\text{Cay}(\mathbb{Z}_p^{2d},\alpha(T)\cup \alpha(T)^{-1})$. Since $\alpha(T)\cup \alpha(T)^{-1}$ is precisely the image of $S_\pm$ under the quotient of $p_{\pm}^{1+2d}$ by $Z(p_\pm^{1+2d})$, the Cayley graphs $\text{Cay}(p_+^{1+2d},S_+)$ and $\text{Cay}(p_-^{1+2d},S_-)$ are covers of $C_p^{2d}$ whose fibers are the cosets of $Z$ in $p_{\pm}^{1+2d}$.

Now suppose there is some 4-cycle $C$ in $\text{Cay}(p_\pm^{1+2d},S_\pm)$. From the definition of a covering graph, the vertices of this 4-cycle must lie in four different fibers, hence the image $\gamma(C)$ under the covering map is also a 4-cycle. The vertex set of any 4-cycle in $\text{Cay}(\mathbb{Z}_p^{2d},\alpha(T)\cup \alpha(T)^{-1})$ is of the form \[\{v,g+v,h+g+v,-g+h+g+v\}\] for some $v\in \mathbb{Z}_p^{2d}$, $g,h \in \alpha(T)$ with $g\notin \{h,h^{-1}\}$. This implies that the vertex set of $C$ is \[\{\tilde{v},\tilde{g}\tilde{v},\tilde{h}\tilde{g}\tilde{v},\tilde{g}^{-1}\tilde{h}\tilde{g}\tilde{v}\}\] for some $\tilde{v}\in\gamma^{-1}(v)$, $\tilde{g},\tilde{h}\in S_\pm$ with $\tilde{g}\notin \{\tilde{h},\tilde{h}^{-1}\}$. Since $C$ is a 4-cycle \[\tilde{v}=\tilde{h}^{-1}\tilde{g}^{-1}\tilde{h}\tilde{g}\tilde{v},\] so $\tilde{g}$ and $\tilde{h}$ commute, contradicting Lemma 5.1 b). 

To conclude, we verify that the two covers are not isomorphic. If a cover of $C_p^{2d}$ contains a $p$-cycle, the image under the covering map of this $p$-cycle must be a cycle of order divisible by $p$, thus a $p$-cycle. The only such $p$-cycles in $\text{Cay}(C_p^{2d},T\cup T^{-1})$ have as vertices the powers of an element of $T$. However, with respect to the multiplication (\ref{minusmult}) each element $t\in S_-$ satisfies $t^p=(0,\dots,0,1)$ if it came from the set $A$, or $t^p=(0,\dots,0,2)$ if it came from the set $B$. So $S_-$ consists only of elements of order $p^2$. It follows that $\text{Cay}(p_-^{1+2d},S_-)$ contains no $p$-cycles.

On the other hand, since each $s\in S_+$ has order $p$, $\text{Cay}(p_+^{1+2d},S_+)$ does contain contains $p$-cycles with vertices \[g,sg,s^2g,\dots, s^{p-1}g\] for each $s\in S_+$ and $g\in G$. So the two covers are not isomorphic.} \qed 

Now it is easy to construct covers for $C_p^{2d-1}$ as induced subgraphs.

\begin{corollary}
Consider $C_p^{2d}$ as the image of the covering maps associated with the two covers given in the previous theorem. If $V$ is the vertex set of an induced $C_p^{2d-1}$ in $C_p^{2d}$, then the preimages of $V$ under the covering maps are the vertex sets of induced non-isomorphic 4-cycle free covers of $C_p^{2d-1}$. \qed
\end{corollary} 

\section{From covers to gain graphs}
In this short section we explicitly describe the the $p$th-root-of-unity valued adjacency matrices for $C_p^{2d}$ that are associated with our covers. This can be done, without too much contrivance, by using the 2-cocycles defined in Lemmas 4.1 and 4.2.

The correspondence between 2-fold coverings of a graph $X=(V,E)$ and $\pm 1$ signings of $X$ can be naturally extended to a correspondence between $r$-fold coverings of $X$ and labelings of $E(X)$ which take values in any permutation group $G$. Define a \textbf{symmetric arc function} to be a function $f: V\times V \rightarrow G\cup\{0\}$ for which $f(u,v)=f(v,u)^{-1}$ whenever $f(u,v)\in G$ and $f(u,v)=0$ if and only if $\{u,v\}\notin E(X)$. The pair $(X,f)$ is called a \textbf{gain graph}, and can be used to define a cover of $X$ as follows. Let $R=\{1,2,\dots, r\}$ be a set on which $G$ acts. The \textbf{cover} $X^f$ of $X$ is the graph with vertex set $V\times R$ and $(v,j)\sim_{Y^f}(u,k)$ exactly when $v\sim_Xu$ and $f(u,v)j=k$. We call $X^f$ the cover \textbf{associated} with $(X,f)$.

\begin{prop}
Let $S$ be the basis for $\mathbb{Z}_p^{2d}$ defined in Section 4, and consider $C_p^{2d}$ as the Cayley graph $\text{Cay}(\mathbb{Z}_p^{2d},S\cup S^{-1})$. Let $\kappa_\pm$ be the 2-cocycles defined in Lemmas 4.1 and 4.2.
For each $s\in S$ and each directed edge $(g,sg)$ in $C_p^{2d}$, define maps $f_\pm:V\times V\rightarrow \mathbb{Z}_p\times \{0\}$ by setting \[f_\pm(g,s+g)=\kappa_\pm(s,g),~~~f_\pm(s+g,g)=\kappa_\pm(s,g)^{-1},\] and then extending so that $f_\pm$ is 0 on all non-arcs of $C_p^d$. Then the gain graphs $(C_p^{2d},f_\pm)$ are associated with covers isomorphic to  $\text{Cay}(p_\pm^{1+2d},S_\pm)$. \qed
\end{prop}

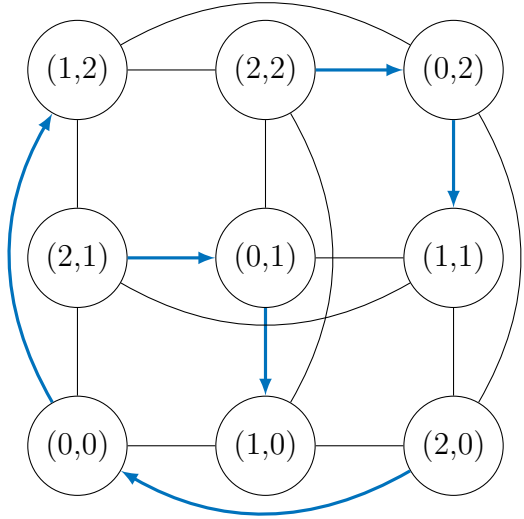
\begin{figure}
\tikzset{>=latex}
    \centering
\begin{tikzpicture}[darkstyle/.style={circle,draw,fill=white!40,minimum size=20}]
  \foreach \x in {0,1,2}
    \foreach \y in {0,1,2} 
       {\pgfmathtruncatemacro{\xlabel}{mod(3+\x-\y,3)}
       \pgfmathtruncatemacro{\ylabel}{\y}
       \node [darkstyle]  (\x\y) at (2.5*\x,2.5*\y) {(\xlabel,\ylabel)};} 

\draw (00)--(10);
\draw (01)[->,very thick, RoyalBlue]--(11);
\draw (02)--(12);

\draw (10)--(20);
\draw (11)--(21);
\draw (12)[->,very thick, RoyalBlue]--(22);

\draw (00)--(01);
\draw [->,very thick, RoyalBlue] (11)--(10);
\draw (20)--(21);

\draw (01)--(02);
\draw (11)--(12);
\draw[->,very thick, RoyalBlue] (22)--(21);

\draw (20) [->,very thick, RoyalBlue,bend left] to (00);
\draw (01) to [bend right] (21);
\draw (02) to [bend left] (22);

\draw (00)[->,very thick, RoyalBlue, bend left] to (02);
\draw (10) to [bend right] (12);
\draw (20) to [bend right] (22);

\end{tikzpicture}
\caption{The gain graph $(C_3^2,f_-).$ Let $\omega$ be a fixed generator of $\mathbb{Z}_3$. The blue arcs are assigned the value $\omega$. The black arcs, drawn as edges, are assigned the value 1. Note that the product of labels around any directed 4-cycle (or directed 3-cycle) is not 1.}

\end{figure}{}

\section{Remarks and loose ends}
We have constructed a new family of graphs that generalize an interesting and useful family. In this final section we consider some of the ways that our new constructions may themselves be interesting or useful. We begin by highlighting some properties of the small examples of our covers, gleaned by computing their automorphism groups. We then consider the prospect of using the spectra of the gain graphs associated with our covers to derive bounds on the maximum degree of induced subgraphs of Cartesian products of cycles.

The first property of the Cohen-Tits cover which we set out to generalize was the presence of its adjacency matrix inside a notable association scheme. See \cite{BCN} Chapter 2 for definitions, and Chapter 9 for the example in question. We had hoped that generalizing the Cohen-Tits covers would yield more association schemes. As of yet, we have found exactly one: It contains the adjacency matrix of our cover $\text{Cay}(3_+^{1+2},S_+)$. This is a symmetric association scheme of rank $7$. The adjacency matrices of our Cayley graphs  $\text{Cay}(3_-^{1+2},S_-)$, $\text{Cay}(3_+^{1+4},S_+)$, $\text{Cay}(3_-^{1+4},S_-)$ sit inside coherent configurations which could perhaps be fused into interesting association schemes. We generated these association schemes and coherent configurations by computing the orbits of a vertex stabilizer of the full automorphism group of the cover. Consequently, the presence of these structures is closely linked to a sufficient amount of symmetry in the covering graphs. On the other hand, we were surprised to find that the covers $\text{Cay}(5_-^{1+2},S_-)$, $\text{Cay}(7_-^{1+2},S_-)$, $\text{Cay}(11_-^{1+2},S_-)$, $\text{Cay}(13_-^{1+2},S_-)$ are as asymmetric as Cayley graphs can possibly be: Their full automorphism groups are respectively $5_-^{1+2}$, $7_-^{1+2}$, $11_-^{1+2}$, and $13_-^{1+2}$. Certainly there is more to say about all of these graphs and the combinatorial structures associated with them.

Given the context of our construction, it is very natural to hope that the existence of these covering graphs allows us to derive meaningful bounds on the maximum degree of an induced subgraph of $C_p^d$. The empirical evidence is not promising.

Let $A$ be the $p$th-root-of-unity valued adjacency matrix for $C_p^{2d}$ associated with the cover $\text{Cay}(p_\pm^{1+2d},S_\pm)$.  The major obstacle is that the size of the largest eigenvalue of $A$ appears to grow very slowly as a function of $d$. Moreover it has relatively low multiplicity. This diminishes the quality of any bound on maximum degree, as well as the size of subgraph to which such a bound can be applied. For instance: Using Huang's method and our covers, we can show that every induced subgraph of $C_3^2$ on at least 4 vertices contains a vertex of degree at least 3, and any induced subgraph of $C_3^4$ on at least 46 vertices has a vertex of degree at least 4. On the one hand, at least for these small examples, our bounds are better than those obtained by Tikaradze \cite{AT}, who considered induced subgraphs of Cartesian products of directed cycles. On the other hand, $C_3^4$ is 3-colorable, so the most interesting problem would appear to be obtaining a bound on the maximum degree of an induced subgraph on $81/3+1=28$ vertices, far outside the reach of our graphs' spectra.

If nothing else, this serves as a reminder of just how perfectly Huang's proof fits its intended purpose. Huang's matrix, or the Cohen-Tits cover, has exactly the correct spectrum to give exactly the right bound on induced subgraphs of exactly the pertinent size. It appears that this is too much to hope for from our generalization.

\section*{Acknowledgement}
I would like to thank Chris Godsil for countless discussions of immeasurable value, and in particular, for suggesting that I look into the Cohen-Tits covers.

Moreover, I would like to a acknowledge partial support from Chris Godsil's research grant NSERC (Canada), Grant No. RGPIN-9439.
\bibliographystyle{plain}

\bibliography{Covers}

\end{document}